\documentclass[10pt,a4paper]{article}

\usepackage{amsfonts,amsthm,amsmath}
\usepackage{graphicx}
\usepackage[height=23cm,width=17cm,twoside,centering]{geometry}
\usepackage{longtable}
\usepackage{array}
\usepackage[arrow,matrix]{xy}

\DeclareMathOperator{\End}{End}

\DeclareMathOperator{\Hom}{Hom}
\DeclareMathOperator{\per}{per}
\DeclareMathOperator{\GL}{GL}

\newcommand{\cC}{\mathcal{C}}
\newcommand{\cD}{\mathcal{D}}
\newcommand{\gL}{\Lambda}
\newcommand{\bQ}{\mathbb{Q}}
\newcommand{\bZ}{\mathbb{Z}}

\theoremstyle{plain}
\newtheorem{thm}{Theorem}[section]
\newtheorem*{thm*}{Theorem}
\newtheorem{prop}[thm]{Proposition}
\newtheorem{lemma}[thm]{Lemma}
\newtheorem{cor}{Corollary}

\theoremstyle{definition}
\newtheorem{dfn}[thm]{Definition}
\newtheorem{remark}[thm]{Remark}
\newtheorem{eg}[thm]{Example}
\newtheorem{notat}[thm]{Notation}
\newtheorem*{quest*}{Question}
\newtheorem{step}{Step}

\begin{document}


\begin{center}
{\Large\sc Derived equivalence classification of the cluster-tilted
algebras of Dynkin type $E$}
\end{center}

\medskip

\begin{center}
{\large\sc Janine Bastian$^*$, Thorsten Holm$^*$\footnote{%
This work has been carried out in the framework of the priority program
SPP 1388 \emph{Darstellungstheorie} of the Deutsche
Forschungsgemeinschaft (DFG). J.\,Bastian and T.\,Holm gratefully
acknowledge financial support through the grant HO 1880/4-1. }, and
Sefi Ladkani$^{**}$\footnote{S.\,Ladkani was supported by a European
Postdoctoral Institute (EPDI) fellowship. He also acknowledges support
from DFG grant LA 2732/1-1 in the framework of the priority program SPP
1388 "Representation theory".} }
\end{center}

\begin{center}
$^*$ Leibniz Universit\"at Hannover, Institut f\"ur Algebra,
Zahlentheorie und Diskrete Mathematik,
Welfengarten 1, 30167 Hannover, Germany\\
{\tt \symbol{`\{}bastian,holm\symbol{`\}}@math.uni-hannover.de}
\end{center}
\begin{center}
$^{**}$ Max-Planck-Institut f\"ur Mathematik, Vivatsgasse 7,
53111 Bonn, Germany\\
Current address: Mathematical Institute of the University of Bonn,
Endenicher Allee 60,\\
53115 Bonn, Germany\\
\tt{sefil@math.uni-bonn.de}
\end{center}


\begin{abstract}
We obtain a complete derived equivalence classification of the
cluster-tilted algebras of Dynkin type $E$. There are 67, 416, 1574
algebras in types $E_6$, $E_7$ and $E_8$ which turn out to fall into 6,
14, 15 derived equivalence classes, respectively. This classification
can be achieved computationally and we outline an algorithm which has been
implemented to carry out this task. We also make the classification
explicit by giving standard forms for each derived equivalence class
as well as complete lists of the algebras contained in each
class; as these lists are quite long they are
provided as supplementary material to this paper.
From a structural point of view the remarkable outcome of our classification
is that two cluster-tilted algebras of Dynkin type $E$ are derived
equivalent if and only if their Cartan matrices represent equivalent
bilinear forms over the integers which in turn happens if and only if
the two algebras are connected by a sequence of ``good'' mutations.
This is reminiscent of the derived equivalence classification of cluster-tilted
algebras of Dynkin type $A$, but quite different from the situation in
Dynkin type $D$ where a far-reaching classification has been obtained
using similar methods as in the present paper but some very subtle
questions are still open.
\end{abstract}


\section{Introduction}

\subsection{The problem}

Cluster categories have been introduced in~\cite{BMRRT} (see
also~\cite{CCS06a} for Dynkin type $A$) as a representation-theoretic
approach to Fomin and Zelevinsky's cluster algebras without
coefficients having skew-symmetric exchange matrices (so that matrix
mutation becomes the combinatorial recipe of mutation of quivers). This
highly successful approach allows to use deep algebraic and
representation-theoretic methods in the context of cluster algebras. A
crucial role is played by the so-called cluster tilting objects in the
cluster category which model the clusters in the cluster algebra. The
endomorphism algebras of these cluster tilting objects are called
cluster-tilted algebras.

Cluster-tilted algebras are particularly well-understood if the quiver
underlying the cluster algebra, and hence the cluster category, is of
Dynkin type. Cluster-tilted algebras of Dynkin type can be described as
quivers with relations where the possible quivers are precisely the
quivers in the mutation class of the Dynkin quiver, and the relations
are uniquely determined by the quiver in an explicit
way~\cite{BMR_finite}. By a result of Fomin and
Zelevinsky~\cite{FominZelevinsky03}, the mutation class of a Dynkin
quiver is finite. Moreover, the quivers in the mutation classes of
Dynkin quivers are explicitly known; for type $A$ they can be found
in~\cite{Buan-Vatne}, for type $D$ in~\cite{Vatne} and for type $E$
they can be enumerated using a computer, for example by the Java
applet~\cite{Keller-software}.

However, despite knowing the cluster-tilted algebras of Dynkin type as
quivers with relations, many structural properties are not understood
yet. One important structural aspect is to understand the derived
module categories of the cluster-tilted algebras. In particular, one
would want to know when two cluster-tilted algebras have equivalent
derived categories. A complete derived equivalence classification has
been achieved so far for cluster-tilted algebras of Dynkin type $A$ by
Buan and Vatne~\cite{Buan-Vatne} and for those of extended Dynkin type
$\tilde{A}$ by the first-named author~\cite{Bastian}. For Dynkin type
$D$, a far-reaching classification has been presented by the authors
in~\cite{BHL10}.
In the present paper we address this problem for cluster-tilted
algebras of Dynkin type $E$ and obtain a complete derived equivalence
classification of these algebras.

\subsection{Main results}

There are two natural approaches to address derived equivalence
classification problems of a given collection of algebras arising from
some combinatorial data. The top-to-bottom approach is to divide the
algebras into equivalence classes according to some invariants of
derived equivalence, so that algebras belonging to different classes
are not derived equivalent. The bottom-to-top approach is to
systematically construct, based on the combinatorial data, derived
equivalences between pairs of these algebras and then to arrange these
algebras into groups where any two algebras are related by a sequence
of such derived equivalences. To obtain a complete derived equivalence
classification one has to combine these approaches and hope that the
two resulting partitions of the entire collection of algebras coincide.

The invariant of derived equivalence we use in this paper is the
integer equivalence class of the bilinear form represented by the
Cartan matrix of an algebra $A$. As this invariant is sometimes
arithmetically subtle to compute directly, we instead compute the
determinant of the Cartan matrix $C_A$ and the characteristic
polynomial of its asymmetry matrix $S_A = C_A C_A^{-T}$, defined
whenever $C_A$ is invertible over $\bQ$, and encode them conveniently
in a single polynomial that we call the \emph{associated polynomial} of
$A$. This quantity is generally a weaker invariant of derived
equivalence, but in our case it will turn out to be enough for the
classification.

The constructions we use are the following. Since any two quivers in a
mutation class are connected by a sequence of mutations, it is natural
to ask when a single mutation of quivers is accompanied by derived
equivalence of their corresponding cluster-tilted algebras. The third
author~\cite{Ladkani10} has presented a procedure to determine when two
cluster-tilted algebras whose quivers are related by a single mutation
are also related by Brenner-Butler (co-)tilting, which is a particular
kind of derived equivalence. We call such quiver mutation \emph{good
mutation}. In other words, a mutation at some vertex is good if the
corresponding Brenner-Butler tilting module is defined and its
endomorphism algebra is isomorphic to the cluster-tilted algebra of the
mutated quiver. Obviously, the cluster-tilted algebras of quivers
connected by a sequence of good mutations are derived equivalent. The
explicit knowledge of the relations for cluster-tilted algebras of
Dynkin type together with the procedure in~\cite{Ladkani10} imply that
for these algebras there is an algorithm to decide if a mutation is
good or not.

It turns out that for cluster-tilted algebras of Dynkin type $E$ the
two approaches can be successfully combined to give a complete derived
equivalence classification.

\begin{thm*}
The following conditions are equivalent for two cluster-tilted algebras
$\gL$ and $\gL'$ of Dynkin type $E$:
\begin{enumerate}
\renewcommand{\theenumi}{\alph{enumi}}
\renewcommand{\labelenumi}{\emph{(\theenumi)}}

\item
$\gL$ and $\gL'$ have the same associated polynomial;

\item
The Cartan matrices of $\gL$ and $\gL'$ represent equivalent bilinear
forms over $\bZ$;

\item
$\gL$ and $\gL'$ are derived equivalent;

\item
The quivers of $\gL$ and $\gL'$ are connected by a sequence of good
mutations.
\end{enumerate}
\end{thm*}

Note that the implication (c) $\Rightarrow $ (b) holds in general for
any two (finite-dimensional) algebras $\gL$ and $\gL'$, and the
implication (b) $\Rightarrow$ (a) holds whenever the associated
polynomials are defined, i.e.\ when the Cartan matrices are invertible
over $\mathbb{Q}$. Moreover, for cluster-tilted algebras the
implication (d) $\Rightarrow $ (c) is evident from the definition.

Let us mention a few consequences from the derived equivalence
classification. Recall that a good mutation involves a particular kind
of derived equivalence, namely Brenner-Butler tilting. However, it
turns out that whenever two cluster-tilted algebras of Dynkin type $E$
whose quivers are related by a single mutation are derived equivalent,
they are also related by Brenner-Butler tilting. In other words:

\begin{cor} \label{c:goodmutder}
A single mutation of quivers in the mutation class of Dynkin type $E$
is good if and only if the corresponding cluster-tilted algebras are
derived equivalent.
\end{cor}

In view of this corollary, condition (d) in the theorem can thus be
rephrased as follows:
\begin{enumerate}
\item [(d')]
\emph{The quivers of $\gL$ and $\gL'$ are connected by a sequence of
mutations such that all the intermediate cluster-tilted algebras are
derived equivalent.}
\end{enumerate}

Another consequence of the classification is the following.
\begin{cor} \label{c:opder}
Let $\gL$ be a cluster-tilted algebra of Dynkin type $E$. Then $\gL$
and its opposite algebra $\gL^{op}$ are derived equivalent.
\end{cor}

In addition to the above general statement we make the derived
equivalence classification explicit by giving standard forms for each
derived equivalence class as well as providing complete lists of the
algebras contained in each class. As these lists are quite long, we
will refrain from reproducing them here. Rather, they can be found in
the supplementary material to this paper, which is freely and
permanently accessible on the arXiv as
part of the text in earlier versions of \cite{BHLTypeE}. It can also be
accessed via the third author's personal web page \cite{Sefi:www}, in a
version generated by using a computer program.

In the following tables we list the associated polynomials of the
cluster-tilted algebras, and also the total number of algebras in each
derived equivalence class.

For type $E_6$ the mutation class consists of 67 quivers. The
corresponding cluster-tilted algebras turn out to fall into six derived
equivalence classes as follows.

\begin{center}
\begin{tabular}{|l|r||l|r|}
\hline
\multicolumn{4}{|c|}{Derived equivalence classes for type $E_6$} \\
\hline
\multicolumn{1}{|c|}{Associated polynomial} & \multicolumn{1}{c||}{\#}  &
\multicolumn{1}{c|}{Associated polynomial} & \multicolumn{1}{c|}{\#} \\
\hline
& & &\\[-10pt]
$x^6-x^5+x^3-x+1$ & 20 & $3(x^6+x^3+1)$ & 19\\
\hline
& & &\\[-10pt]
$2(x^6-x^4+2x^3-x^2+1)$ & 16 & $4(x^6+x^4+x^2+1)$ & 7 \\
\hline
& & &\\[-10pt]
$2(x^6-2x^4+4x^3-2x^2+1)$ & 3 & $4(x^6+x^5-x^4+2x^3-x^2+x+1)$ & 2 \\
\hline
\end{tabular}
\end{center}

For type $E_7$ the mutation class consists of 416 quivers. The derived
equivalence classes of the corresponding cluster-tilted algebras are
again characterized by the associated polynomials; there are 14 classes
in total, given as follows.
\begin{center}
\begin{tabular}{|l|r||l|r|}
\hline
\multicolumn{4}{|c|}{Derived equivalence classes for type $E_7$} \\
\hline
\multicolumn{1}{|c|}{Associated polynomial} & \multicolumn{1}{c||}{\#}  &
\multicolumn{1}{c|}{Associated polynomial} & \multicolumn{1}{c|}{\#} \\
\hline
& & &\\[-10pt]
$x^7-x^6+x^4-x^3+x-1$ & 64 & $4(x^7+x^6-2x^5+2x^4-2x^3+2x^2-x-1)$ & 2 \\
\hline
& & &\\[-10pt]
$2(x^7-x^5+2x^4-2x^3+x^2-1)$ & 32 & $4(x^7+x^5-x^4+x^3-x^2-1)$ & 56 \\
\hline
& & &\\[-10pt]
$2(x^7-x^5+x^4-x^3+x^2-1)$ & 72 & $4(x^7+x^5-2x^4+2x^3-x^2-1)$ & 8 \\
\hline
& & &\\[-10pt]
$2(x^7-2x^5+4x^4-4x^3+2x^2-1)$ & 8 & $5(x^7+x^5-x^4+x^3-x^2-1)$ & 17 \\
\hline
& & &\\[-10pt]
$3(x^7-1)$ & 124 & $6(x^7+x^6-x^4+x^3-x-1)$ & 11 \\
\hline
& & &\\[-10pt]
$4(x^7+x^6-x^5+x^4-x^3+x^2-x-1)$ & 16 & $6(x^7+x^5-x^2-1)$ & 1 \\
\hline
& & &\\[-10pt]
$4(x^7+x^6-x^5-x^4+x^3+x^2-x-1)$ & 4 & $8(x^7+x^6+x^5-x^4+x^3-x^2-x-1)$ & 1   \\
\hline
\end{tabular}
\end{center}

For type $E_8$ the mutation class consists of 1574 quivers. The
corresponding cluster-tilted algebras turn out to fall into 15
different derived equivalence classes which are characterized as
follows.
\begin{center}
\begin{tabular}{|l|r||l|r|}
\hline
\multicolumn{4}{|c|}{Derived equivalence classes for type $E_8$} \\
\hline
\multicolumn{1}{|c|}{Associated polynomial} & \multicolumn{1}{c||}{\#}  &
\multicolumn{1}{c|}{Associated polynomial} & \multicolumn{1}{c|}{\#} \\
\hline
& & &\\[-10pt]
$x^8-x^7+x^5-x^4+x^3-x+1$ & 128 & $4(x^8+x^6-x^5+2x^4-x^3+x^2+1)$ & 221 \\
\hline
& & &\\[-10pt]
$2(x^8-x^6+2x^5-2x^4+2x^3-x^2+1)$ & 64 & $4(x^8+x^6-2x^5+4x^4-2x^3+x^2+1)$ & 22 \\
\hline
& & &\\[-10pt]
$2(x^8-x^6+x^5+x^3-x^2+1)$ & 256 & $5(x^8+x^6+x^4+x^2+1)$ & 167 \\
\hline
& & &\\[-10pt]
$2(x^8-2x^6+4x^5-4x^4+4x^3-2x^2+1)$ & 16 & $6(x^8+x^6+x^5+x^3+x^2+1)$ & 38 \\
\hline
& & &\\[-10pt]
$3(x^8+x^4+1)$ & 384 & $6(x^8+x^7+2x^4+x+1)$ & 118 \\
\hline
& & &\\[-10pt]
$4(x^8+x^7-x^6+x^5+x^3-x^2+x+1)$ & 72 & $8(x^8+2x^7+2x^4+2x+1)$ & 4 \\
\hline
& & &\\[-10pt]
$4(x^8+x^7-x^6+2x^4-x^2+x+1)$ & 48 & $8(x^8+x^7+x^6+2x^4+x^2+x+1)$ & 24   \\
\hline
& & &\\[-10pt]
$4(x^8+x^7-2x^6+2x^5+2x^3-2x^2+x+1)$ & 12 & & \\
\hline
\end{tabular}
\end{center}

\subsection{Comparison with the other Dynkin types}

We now put our results in perspective by comparing them to the derived
equivalence classifications of the cluster-tilted algebras of the other
Dynkin types. In type $A$, a complete classification has been achieved
by Buan and Vatne~\cite{Buan-Vatne}, whereas in type $D$, a
far-reaching classification has been presented by the authors
in~\cite{BHL10}.

It turns out that two cluster-tilted algebras of type $A_n$ are derived
equivalent if and only if their quivers have the same number of
3-cycles. For distinguishing such algebras up to derived equivalence
one uses the determinants of the Cartan matrices; these have been
determined explicitly for arbitrary gentle algebras by the second
author in~\cite{Holm-gentle}.

A quick look at the above tables reveals that in Dynkin type $E$,
unlike as in type $A$, the determinant itself is not sufficient for
distinguishing the algebras up to derived equivalence. However, it is
interesting to note that two cluster-tilted algebras of type $E_n$
($n=6,7,8$) with the same \emph{odd} Cartan determinant are derived
equivalent.

For this reason one needs to use further derived invariants, such as
the characteristic polynomial of the asymmetry matrix of the Cartan
matrix, in order to distinguish the derived equivalence classes.
Looking again at the above tables we see that this polynomial alone is
not enough to distinguish the derived classes: there are two derived
equivalence classes of cluster-tilted algebras of type $E_7$ with the
same characteristic polynomial $x^7+x^5-x^4+x^3-x^2-1$. Our main
theorem claims that by combining this polynomial with the Cartan
determinant to form the associated polynomial we get a complete
invariant of derived equivalence for cluster-tilted algebras of type
$E$.

There are five non-trivial assertions which are satisfied by the
cluster-tilted algebras of type $E$. These include the three
implications (a) $\Rightarrow$ (b), (b) $\Rightarrow$ (c) and (c)
$\Rightarrow$ (d) in the main theorem, together with the statements of
Corollaries~\ref{c:goodmutder} and~\ref{c:opder}. In the table below we
specify, for each of these assertions and each of the Dynkin types $A$,
$D$ and $E$, whether the assertion holds ($\surd$), does not hold
($\times$) or is unknown (?) for cluster-tilted algebras of the given
type. In the two latter cases, we also give the smallest numbers of
vertices for which there is a counterexample (in the case '$\times$')
or an example that could not be settled (in the case '?').

\begin{center}
\begin{tabular}{|c||c|cc|c||}
\multicolumn{1}{c||}{} &
Type $A$ & \multicolumn{2}{c|}{Type $D$} & Type $E$ \\ \hline
(a) $\Rightarrow$ (b) &
$\surd$ & $\times$ & {\footnotesize ($D_{15}$)} & $\surd$ \\ \hline
(b) $\Rightarrow$ (c) &
$\surd$ & ? & {\footnotesize ($D_{15}, D_{19}$)} & $\surd$ \\ \hline
(c) $\Rightarrow$ (d) &
$\surd$ & $\times$ & {\footnotesize ($D_6, D_8$)} & $\surd$ \\ \hline
Corollary~\ref{c:goodmutder} &
$\surd$ & \multicolumn{2}{l|}{$\surd$} & $\surd$ \\ \hline
Corollary~\ref{c:opder} &
$\surd$ & ? & {\footnotesize ($D_{15}$)} & $\surd$ \\ \hline
\end{tabular}
\end{center}

We see that as far as these assertions are concerned, cluster-tilted
algebras of Dynkin type $E$ behave similarly to that of type $A$. For
type $D$, however, the picture is quite different; already in types
$D_6$ and $D_8$ there are pairs of derived equivalent cluster-tilted
algebras whose quivers cannot be connected by sequences of good
mutations, making it necessary to use further constructions of derived
equivalence. Thus, the derived equivalence classification of
cluster-tilted algebras of type $E$ is simpler than that of type $D$,
even for small numbers of vertices.

Furthermore, starting at type $D_{15}$ there are pairs of
cluster-tilted algebras which we could not decide whether they are
derived equivalent, or not; there does not seem to be an accessible
tilting complex for these algebras but on the other hand all computable
derived invariants available to us are the same for both algebras. It
is therefore a remarkable coincidence that the two effectively
decidable conditions (a) and~(d) are equivalent for cluster-tilted
algebras of type $E$, thus enabling a complete derived equivalence
classification of these algebras.

\subsection*{Acknowledgement}

The first two authors are grateful to Bernhard Keller for suggesting to
contact the third author on questions left open in a first version of
this paper which could then be completely settled in the present joint
version.

\section{Assembling the proof}

Throughout this paper let $K$ be an algebraically closed field. All
algebras are assumed to be finite-dimensional $K$-algebras. For an
algebra $A$, we denote the bounded derived category of right
$A$-modules by $\cD^b(A)$. Two algebras $A$ and $B$ are called
\emph{derived equivalent} if $\cD^b(A)$ and $\cD^b(B)$ are equivalent
as triangulated categories.

\subsection{The equivalence class of the Euler form as derived invariant}

Let $A$ be an algebra and let $P_1, \dots, P_n$ be a complete
collection of non-isomorphic indecomposable projective right
$A$-modules (finite-dimensional over $K$). The \emph{Cartan matrix} of
$A$ is then the $n \times n$ matrix $C_A$ defined by $(C_A)_{ij} =
\dim_K \Hom_A(P_j, P_i)$.

Denote by $\per A$ the triangulated category of \emph{perfect}
complexes of $A$-modules inside the derived category of $A$, that is,
the complexes which are quasi-isomorphic to finite complexes of
finitely generated projective $A$-modules. The Grothendieck group
$K_0(\per A)$ is a free abelian group on the generators $[P_1], \dots,
[P_n]$, and the expression
\[
\langle X, Y \rangle = \sum_{r \in \bZ} (-1)^r \dim_K \Hom_{\per A}(X,
Y[r])
\]
is well defined for any $X, Y \in \per A$ and induces a bilinear form
on $K_0(\per A)$, known as the \emph{Euler form}, whose matrix with
respect to the basis of projectives is $C_A^T$.

The following proposition is well known. For the convenience of the
reader, we give the short proof, see also the \emph{proof} of
Proposition~1.5 in~\cite{Bocian}.

\begin{prop}
Let $A$ and $B$ be two finite-dimensional, derived equivalent algebras.
Then the matrices $C_A$ and $C_B$ represent equivalent bilinear forms
over $\bZ$, that is, there exists $P \in \GL_n(\bZ)$ such that $P C_A
P^T = C_B$, where $n$ denotes the number of non-isomorphic
indecomposable projective modules of $A$ (and $B$).
\end{prop}
\begin{proof}
Indeed, by~\cite{Rickard}, if $A$ and $B$ are derived equivalent, then
$\per A$ and $\per B$ are equivalent as triangulated categories. Now
any triangulated functor $F: \per A \to \per B$ induces a linear map
from $K_0(\per A)$ to $K_0(\per B)$. When $F$ is also an equivalence,
this map is an isomorphism of the Grothendieck groups preserving the
Euler forms. Thus, if $[F]$ denotes the matrix of this map with respect
to the bases of indecomposable projectives, then $[F]^T C_B [F] = C_A$.
\end{proof}

In general, to decide whether two integral bilinear forms are
equivalent is a very subtle arithmetical problem. Therefore, it is
useful to introduce somewhat weaker invariants that are computationally
easier to handle. In order to do this, assume further that $C_A$ is
invertible over $\bQ$, that is, $\det C_A \neq 0$. In this case one can
consider the rational matrix $S_A = C_A C_A^{-T}$ (here $C_A^{-T}$
denotes the inverse of the transpose of $C_A$), known in the theory of
non-symmetric bilinear forms as the \emph{asymmetry} of $C_A$.

\begin{prop}
Let $A$ and $B$ be two finite-dimensional, derived equivalent algebras
with invertible (over $\bQ$) Cartan matrices. Then we have the
following assertions, each implied by the preceding one:
\begin{enumerate}
\renewcommand{\theenumi}{\alph{enumi}}
\renewcommand{\labelenumi}{\emph{(\theenumi)}}
\item
There exists $P \in \GL_n(\bZ)$ such that $P C_A P^T = C_B$.

\item
There exists $P \in \GL_n(\bZ)$ such that $P S_A P^{-1} = S_B$.

\item
There exists $P \in \GL_n(\bQ)$ such that $P S_A P^{-1} = S_B$.

\item
The matrices $S_A$ and $S_B$ have the same characteristic polynomial.
\end{enumerate}
\end{prop}

For proofs and discussion, see for example \cite[Section~3.3]{Ladkani}.
Since the determinant of an integral bilinear form is also invariant
under equivalence, we obtain the following discrete invariant of
derived equivalence.

\begin{dfn}
For an algebra $A$ with invertible Cartan matrix $C_A$ over $\bQ$, we
define its \emph{associated polynomial} as $(\det C_A) \cdot
\chi_{S_A}(x)$, where $\chi_{S_A}(x)$ is the characteristic polynomial
of the asymmetry matrix $S_A = C_A C_A^{-T}$.
\end{dfn}

\begin{remark}
The matrix $S_A = C_A C_A^{-T}$ (or better, minus its transpose
$-C_A^{-1} C_A^T$) is related to the \emph{Coxeter transformation}
which has been widely studied in the case when $A$ has finite global
dimension (so that $C_A$ is invertible over $\bZ$), see~\cite{Lenzing}.
It is then the $K$-theoretic shadow of the Serre functor and the
related Auslander-Reiten translation in the derived category. The
characteristic polynomial is then known as the \emph{Coxeter
polynomial} of the algebra.
\end{remark}

\begin{remark}
In general, $S_A$ might have non-integral entries. However, when the
algebra $A$ is \emph{Gorenstein}, the matrix $S_A$ is integral, which
is an incarnation of the fact that the injective modules have finite
projective resolutions. By a result of Keller and Reiten
\cite{Keller-Reiten}, this is the case for cluster-tilted algebras.
\end{remark}

\subsection{Mutations of algebras}
We recall the notion of mutations of algebras from~\cite{Ladkani10}.
These are local operations on an algebra $A$ producing new algebras
derived equivalent to $A$.

Let $A = KQ/I$ be an algebra given as a quiver with relations. For any
vertex $i$ of $Q$, there is a trivial path $e_i$ of length 0; the
corresponding indecomposable projective $P_i=e_i A$ is spanned by the
images of the paths starting at $i$. Thus an arrow $i
\xrightarrow{\alpha} j$ gives rise to a map $P_j \to P_i$ given by left
multiplication with $\alpha$.

Let $k$ be a vertex of $Q$ without loops. Consider the following two
complexes of projective $A$-modules
\begin{align*}
T^-_k(A) = \bigl( P_k \xrightarrow{f} \bigoplus_{j \to k} P_j \bigr)
           \oplus \bigl( \bigoplus_{i \neq k} P_i \bigr) &, &
T^+_k(A) = \bigl( \bigoplus_{k \to j} P_j \xrightarrow{g} P_k \bigr)
           \oplus \bigl( \bigoplus_{i \neq k} P_i \bigr)
\end{align*}
where the map $f$ is induced by all the maps $P_k \to P_j$
corresponding to the arrows $j \to k$ ending at $k$, the map $g$ is
induced by the maps $P_j \to P_k$ corresponding to the arrows $k \to j$
starting at $k$, the term $P_k$ lies in degree $-1$ in $T^-_k(A)$ and
in degree $1$ in $T^+_k(A)$, and all other terms are in degree $0$.

\begin{dfn}
Let $A$ be an algebra given as a quiver with relations and $k$ a vertex
without loops.
\begin{enumerate}
\renewcommand{\theenumi}{\alph{enumi}}
\renewcommand{\labelenumi}{(\theenumi)}
\item
We say that the negative mutation of $A$ at $k$ is \emph{defined} if
$T^-_k(A)$ is a tilting complex over $A$. In this case, we call
the algebra
$\mu^-_k(A) = \End_{\cD^b(A)} T^-_k(A)$ the \emph{negative mutation} of
$A$ at the vertex $k$.

\item
We say that the positive mutation of $A$ at $k$ is \emph{defined} if
$T^+_k(A)$ is a tilting complex over $A$. In this case, we call
the algebra
$\mu^+_k(A) = \End_{\cD^b(A)} T^+_k(A)$ the \emph{positive mutation} of
$A$ at the vertex $k$.
\end{enumerate}
\end{dfn}

\begin{remark}
By Rickard's Morita theory for derived categories~\cite{Rickard}, the
negative and the positive mutations of an algebra $A$ at a vertex, when
defined, are always derived equivalent to $A$.
\end{remark}

\begin{eg}
When $k$ is a sink in the quiver of $A$, then $\mu^-_k(A)$ is defined
and moreover $T^-_k(A)$ is isomorphic in $\cD^b(A)$ to the APR-tilting
module~\cite{APR79} corresponding to $k$. Similarly, when $k$ is a
source then $\mu^+_k(A)$ is defined.
\end{eg}

There is a combinatorial criterion to determine whether a mutation at a
vertex is defined, see~\cite[Prop.~2.3]{Ladkani10}. Since the algebras
we will be dealing with in this paper are schurian, we state here the
criterion only for this case, as it takes a particularly simple form.
Recall that an algebra is \emph{schurian} if the entries of its Cartan
matrix are only $0$ or $1$.

\begin{prop} \label{p:critmut}
Let $A$ be a schurian algebra and let $k$ be a vertex in the quiver of
$A$.
\begin{enumerate}
\renewcommand{\theenumi}{\alph{enumi}}
\renewcommand{\labelenumi}{\emph{(\theenumi)}}
\item
The negative mutation $\mu^-_k(A)$ is defined if and only if for any
non-zero path $k \rightsquigarrow i$ starting at $k$ and ending at some
vertex $i$, there exists an arrow $j \to k$ such that the composition
$j \to k \rightsquigarrow i$ is non-zero in $A$.

\item
The positive mutation $\mu^+_k(A)$ is defined if and only if for any
non-zero path $i \rightsquigarrow k$ starting at some vertex $i$ and
ending at $k$, there exists an arrow $k \to j$ such that the
composition $i \rightsquigarrow k \to j$ is non-zero in $A$.
\end{enumerate}
\end{prop}

\begin{remark}
It follows from~\cite[Remark~2.10]{Ladkani10} that in many cases, and
in particular when $A$ is schurian, the negative mutation of $A$ at $k$
is defined if and only if one can associate with $k$ the corresponding
Brenner-Butler tilting module~\cite{BrennerButler80}. Moreover,
$T^-_k(A)$ is then isomorphic in $\cD^b(A)$ to that Brenner-Butler
tilting module.
\end{remark}

\begin{remark} \label{rem:mutopp}
It is enough to introduce only one of the notions of negative and
positive mutation, since the other can then be defined using the notion
of the opposite algebra. Indeed, if $Q$ is the quiver of an algebra
$A$, then the quiver of the opposite algebra $A^{op}$ is the opposite
quiver $Q^{op}$ which has the same vertices as $Q$ and is obtained from
it by inverting all the arrows. The positive mutation $\mu^+_k(A)$ is
then defined if and only if the negative mutation $\mu^-_k(A^{op})$ is
defined. Moreover, in that case, we have $\mu^+_k(A) \simeq \left(
\mu^-_k(A^{op})\right)^{op}$.
\end{remark}

\subsection{Cluster-tilted algebras}

In this section we assume that all quivers are without loops and
$2$-cycles. Given such a quiver $Q$ and a vertex $k$, we denote by
$\mu_k(Q)$ the Fomin-Zelevinsky quiver
mutation~\cite{FominZelevinsky02} of $Q$ at $k$. Two quivers are called
\emph{mutation equivalent} if one can be reached from the other by a
finite sequence of quiver mutations. The \emph{mutation class} of a
quiver $Q$ is the set of all quivers which are mutation equivalent to
$Q$.

For a quiver $Q'$ without oriented cycles, the corresponding cluster
category $\cC_{Q'}$ was introduced in~\cite{BMRRT}. A
\emph{cluster-tilted algebra} of \emph{type $Q'$} is an endomorphism
algebra of a cluster-tilting object in $\cC_{Q'}$, see~\cite{BMR}. It
is known by~\cite{BMR} that for any quiver $Q$ mutation equivalent to
$Q'$, there is a cluster-tilted algebra whose quiver is $Q$. Moreover,
by~\cite{BIRSm}, it is unique up to isomorphism. Hence, there is a
bijection between the quivers in the mutation class of an acyclic
quiver $Q'$ and the isomorphism classes of cluster-tilted algebras of
type $Q'$. This justifies the following notation.

\begin{notat}
Throughout the paper, for a quiver $Q$ which is mutation equivalent to
an acyclic quiver, we denote by $\gL_Q$ the corresponding
cluster-tilted algebra.
\end{notat}

When $Q'$ is a Dynkin quiver of types $A$, $D$ or $E$, the
corresponding cluster-tilted algebras are said to be of Dynkin type.
These algebras have been investigated in~\cite{BMR_finite}, where it is
shown that they are schurian and moreover they can be defined by using
only zero and commutativity relations that can be extracted from their
quivers in an algorithmic way.

\subsection{Good quiver mutations}

For simplicity, we assume in this section that all the quivers we deal
with are mutation equivalent to Dynkin quivers. A more general
treatment can be found in~\cite{Ladkani10}.

Let $Q$ be such a quiver and let $k$ be a vertex of $Q$. Starting with
the cluster-tilted algebra $\gL_Q$, there are two notions of mutation
that one may consider. The first is mutation of quivers, leading to the
cluster-tilted algebra $\gL_{\mu_k(Q)}$. The second is mutation of
algebras, leading (when defined) to the algebra $\mu^-_k(\gL_Q)$ or
$\mu^+_k(\gL_Q)$. It is interesting to ask when these two notions are
compatible. This motivates the following definition.

\begin{dfn}
Let $Q$ be a quiver which is mutation equivalent to a Dynkin quiver and
let $k$ be a vertex of $Q$. The mutation of $Q$ at $k$ is \emph{good}
if $\gL_{\mu_k(Q)} \simeq \mu^-_k(\gL_Q)$ or $\gL_{\mu_k(Q)} \simeq
\mu^+_k(\gL_Q)$ (or both).
\end{dfn}

Obviously, a good mutation implies the derived equivalence of the
corresponding cluster-tilted algebras.

\begin{remark}
Observe that the condition $\gL_{\mu_k(Q)} \simeq \mu^+_k(\gL_Q)$ is
equivalent to the condition $\gL_Q \simeq \mu^-_k(\gL_{\mu_k(Q)})$, so
that the mutation of $Q$ at $k$ is good if and only if that of
$\mu_k(Q)$ at $k$ is good. Hence being ``good'' is a property of the
mutation, independently of which of the two quivers we start with.
\end{remark}

\emph{A-priori}, checking a condition like $\gL_{\mu_k(Q)} \simeq
\mu^-_k(\gL_Q)$ involves two steps: the first is to check that the
algebra mutation $\mu^-_k(\gL_Q)$ is defined, and the second is to
check that it is isomorphic to $\gL_{\mu_k(Q)}$. This last step is not
so easily adapted for a computer since it involves the computation of
an endomorphism algebra of a tilting complex.

For cluster-tilted algebras of Dynkin type, the statement of
Theorem~5.3 in~\cite{Ladkani10}, linking more generally mutation of
cluster-tilting objects in 2-Calabi-Yau categories with mutations of
their endomorphism algebras, takes the following form.

\begin{prop} \label{p:goodmut}
Let $Q$ be mutation equivalent to a Dynkin quiver and let $k$ be a
vertex of $Q$.
\begin{enumerate}
\renewcommand{\theenumi}{\alph{enumi}}
\renewcommand{\labelenumi}{\emph{(\theenumi)}}
\item
$\gL_{\mu_k(Q)} \simeq \mu^-_k(\gL_Q)$ if and only if the two algebra
mutations $\mu^-_k(\gL_Q)$ and $\mu^+_k(\gL_{\mu_k(Q)})$ are defined.

\item
$\gL_{\mu_k(Q)} \simeq \mu^+_k(\gL_Q)$ if and only if the two algebra
mutations $\mu^+_k(\gL_Q)$ and $\mu^-_k(\gL_{\mu_k(Q)})$ are defined.
\end{enumerate}
\end{prop}

Thus, in order to check the conditions in the definition of good
mutation it is enough to check that certain algebra mutations are
defined.

\begin{remark}
In view of Propositions~\ref{p:critmut} and~\ref{p:goodmut}, there is
an algorithm which decides, given a quiver which is mutation equivalent
to a Dynkin quiver, whether a mutation at a vertex is good or not.
\end{remark}

We recollect a few basic observations on good mutations.

\begin{lemma} \label{l:goodmutss}
A mutation of $Q$ at a sink or a source is always good.
\end{lemma}
\begin{proof}
Assume that $k$ is a sink in $Q$. Then the negative mutation
$\mu^-_k(\gL_Q)$ is defined. Since $k$ becomes a source in the quiver
$\mu_k(Q)$, the positive mutation $\mu^+_k(\gL_{\mu_k(Q)})$ is defined.
Now apply Proposition~\ref{p:goodmut}.
\end{proof}

Hence, for the purpose of derived equivalence classification it is
enough to consider quivers only up to \emph{sink/source equivalence}:
two quivers are called sink/source equivalent if one can be obtained
from the other by performing mutations only at vertices which are sinks
or sources.

\begin{lemma} \label{l:goodmutop}
The mutation of $Q$ at a vertex $k$ is good if and only if the mutation
of $Q^{op}$ at $k$ is good.
\end{lemma}
\begin{proof}
Observe that $\mu_k(Q^{op}) \simeq \left(\mu_k(Q)\right)^{op}$ and
$\gL_{Q^{op}} \simeq \gL_Q^{op}$. Now use Remark~\ref{rem:mutopp} and
Proposition~\ref{p:goodmut}.
\end{proof}

\subsection{Putting everything together}

We now have all the ingredients needed for the proof of the main
theorem. It is enough to show the equivalence of conditions~(a) and~(d)
in that theorem, i.e.\ that the quivers of any two cluster-tilted
algebras of Dynkin type $E$ with the same associated polynomial can be
connected by a sequence of good mutations.

The proof is computational; we give below an outline of an algorithm to
carry out this task. We start with a Dynkin quiver $Q$ of type $E_n$
where $n \in \{6, 7, 8\}$.

\begin{step}[Mutation data]
Compute the mutation class $Q_1, Q_2, \dots, Q_N$ of $Q$ together with
functions
\begin{align*}
s : \{1, \dots, N\} \times \{1, \dots, n\} \to \{1, \dots, N\} &,& \pi
:\{1, \dots, N\} \times \{1, \dots, n\} \to S_n
\end{align*}
such that for each index $1 \leq t \leq N$ and vertex $1 \leq k \leq n$
we have $\mu_k(Q_t) \simeq Q_{s(t,k)}$ with the isomorphism given by
the permutation $\pi(t,k)$ on the set of vertices. This is done as follows:
\begin{enumerate}
\item
Set $Q_1 \leftarrow Q$, $t \leftarrow 1$, $m \leftarrow 1$.

\item \label{it:tstg}
At the $t$-th stage, we have a list $Q_1, \dots, Q_m$ for some $m \geq
t$. For each $1 \leq k \leq n$, compute $\mu_k(Q_t)$. If it is
isomorphic to some member in the list, set $s(t,k)$ and $\pi(t,k)$
accordingly; otherwise, append it to the list and set $m \leftarrow
m+1$.

\item
If $m=t$, set $N \leftarrow m$, stop; otherwise, set $t \leftarrow t+1$
and go to stage~\ref{it:tstg}.
\end{enumerate}
\end{step}

For $1 \leq t \leq N$, denote by $\gL_t$ the cluster-tilted algebra of
$Q_t$ and by $C_t$ its Cartan matrix. Perform Steps~\ref{st:nzpath},
\ref{st:Cartan} and~\ref{st:tilting} below for each $1 \leq t \leq N$:

\begin{step}[Algebras and non-zero paths] \label{st:nzpath}
Compute the zero-relations and the commutativity-relations of $\gL_t$
from the quiver $Q_t$ according to the algorithm in~\cite{BMR_finite}.
Then, for each $1 \leq i, j \leq n$, store the list $L_t(i,j)$ of paths
starting at $i$ and ending at $j$ whose image in $\gL_t$ is non-zero.
Note that since $\gL_t$ is schurian, in computing these lists it is
enough to consider paths traversing through each vertex at most once.
\end{step}

\begin{step}[Cartan matrix and associated polynomial] \label{st:Cartan}
For each $1 \leq i,j \leq n$, set $(C_t)_{ij} = 1$ if $L_t(i,j)$ is not
empty, otherwise set $(C_t)_{ij} = 0$. Then compute the asymmetry $S_t
= C_t C_t^{-T}$ and the associated polynomial $p_t(x) = (\det C_t)
\cdot \chi_{S_t}(x)$.
\end{step}

\begin{step}[Tilting data] \label{st:tilting}
For each vertex $1 \leq k \leq n$, check which of the algebra mutations
$\mu^-_k(\gL_t)$ and $\mu^+_k(\gL_t)$ is defined; to this end use the
lists $L_t(i,j)$ and apply the criteria of Proposition~\ref{p:critmut}.
\end{step}

We encode the good mutations in an undirected graph $G$ whose vertices
are the indices $1,2,\dots, N$ and for any $1 \leq t,t' \leq N$ there
is an edge $t - t'$ if and only if the quivers $Q_t$ and $Q_{t'}$ are
related by a good mutation at some vertex. The graph $G$ is computed in
the next step as follows:

\begin{step}[Good mutations graph]
For each $1 \leq t \leq N$ and $1 \leq k \leq n$:
\begin{itemize}
\item
Set $t' \leftarrow s(t,k)$ and $\ell \leftarrow \pi(t,k)^{-1}(k)$;

\item
If $\mu^-_k(\gL_t)$ and $\mu^+_{\ell}(\gL_{t'})$ are defined, set an
edge $t - t'$ in $G$;

\item
If $\mu^+_k(\gL_t)$ and $\mu^-_{\ell}(\gL_{t'})$ are defined, set an
edge $t - t'$ in $G$.
\end{itemize}
\end{step}

Observe that two quivers $Q_t$ and $Q_{t'}$ are connected by a sequence
of good mutations if and only if the indices $t$ and $t'$ are in the
same connected component of the graph $G$.

\begin{step}[End of proof]
Compute the connected components of the graph $G$ and then take one
representative from each connected component (e.g.\ the smallest one in
that component) to get a sequence $t_1, t_2, \dots, t_M$. Finally
verify that the associated polynomials $p_{t_1}(x), \dots, p_{t_M}(x)$
are all distinct.
\end{step}

We have made two implementations of this algorithm. The first is a
fully automatic one developed by the third author by writing computer
code for the \textsc{Magma} computational algebra system~\cite{Magma}.
The program takes a Dynkin quiver as input and produces all the
relevant information (quivers in the mutation class, defining relations
for the cluster-tilted algebras, Cartan matrices, lists of good
mutations and the partition according to the associated polynomials).
Its output can be found in the supplementary material to this paper,
available on the third author's personal web page \cite{Sefi:www}.
Since all the mutations are examined by the program, it can also verify
Corollary~\ref{c:goodmutder} in addition to the main theorem.

The other implementation is more ``traditional''. The mutation class is
computed using Keller's Java applet~\cite{Keller-software}, then
Steps~\ref{st:nzpath} and~\ref{st:Cartan} are performed and the quivers
are divided into groups according to the associated polynomials of
their corresponding cluster-tilted algebras. Then one finds enough good
mutations to show that in each such group, any two quivers can be
connected by good mutations.

In this approach one does not build the graph $G$ completely, but
rather finds enough of its edges to construct a spanning tree for each
set of vertices with the same associated polynomial. Shortcuts to
reduce the number of calculations are provided by
Lemma~\ref{l:goodmutss} (sink/source equivalence) and
Lemma~\ref{l:goodmutop} (opposites). We shall demonstrate this approach
for the quiver $E_6$ in the next section. For the quivers $E_7$ and
$E_8$, lists of good mutations can be found in a previous version
of this paper \cite{BHLTypeE}.

\section{Cluster-tilted algebras of type $E_6$}
\label{Sec:E6}

In this section we demonstrate the techniques of the previous sections
and describe in detail the derived equivalence classification of
cluster-tilted algebras of type $E_6$. Since the number of these
algebras is not large, we can present the full classification and give
enough details so that the reader may reproduce all our findings.

We start by computing the mutation class of the $E_6$ quiver given by
\begin{center}
\includegraphics[scale=1.2]{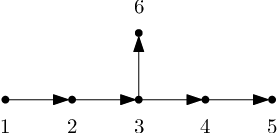}
\end{center}
This can be done, for example, by using the Java applet of
Keller~\cite{Keller-software}. The mutation class of $E_6$ consists of
67 quivers. To reduce the number of calculations, observe that for the
purpose of derived equivalence classification of the corresponding
cluster-tilted algebras it suffices to consider the quivers only up to
sink/source equivalence, see Lemma~\ref{l:goodmutss}. There are $21$
quivers up to sink/source equivalence, and we number them
$1,2,\dots,21$ according to the output of the
program~\cite{Keller-software}.

For each representative quiver in a sink/source equivalence class we
compute the relations of the corresponding cluster-tilted algebra
according to~\cite{BMR_finite}. From this we deduce the Cartan matrix
and the associated polynomial, obtained by multiplying the determinant
of the Cartan matrix by the characteristic polynomial of its asymmetry
matrix.

As there are six distinct such polynomials, the quivers are divided
into six groups according to the associated polynomial of their
corresponding cluster-tilted algebras. We claim that these groups form
the six derived equivalence classes of the cluster-tilted algebras of
type $E_6$.

In the table below we list, for each of these six polynomials, all the
quivers in the mutation class of type $E_6$ whose corresponding
cluster-tilted algebras have this polynomial as associated polynomial.
For each set of sink/source equivalent quivers we give only one picture
where certain arrows are replaced by undirected lines; this has to be
read that these lines can take any orientation. We also give the
relations of the corresponding cluster-tilted algebras. Since there is
at most one arrow between any two vertices, we indicate a path by the
sequence of vertices it traverses. A zero-relation is then indicated by
a sequence of the form $(a,b,c,\dots)$ and a commutativity relation has
the form $(a,b,c,\dots)-(a',b',c',\dots)$.


{\scriptsize
\begin{center}
\begin{tabular}[t]{|c|p{5cm}<{\centering}|p{5cm}<{\centering}|}
\hline
\multicolumn{3}{|c|}{}\\[-6pt]
\multicolumn{3}{|c|}{$\boldsymbol{x^6 - x^5 + x^3 - x + 1}$}\\
\hline
no. & quiver & relations\\
\hline
$1$ & \begin{tabular}{c} \\ \includegraphics[scale=0.8]{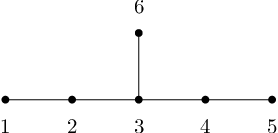} \\ \\ \end{tabular} &
\begin{tabular}{c} \\
None \\ \end{tabular}\\
\hline
\end{tabular}
\end{center}}

{\scriptsize
\begin{center}
\begin{tabular}[t]{|c|p{5cm}<{\centering}|p{5cm}<{\centering}|}
\hline
\multicolumn{3}{|c|}{}\\[-6pt]
\multicolumn{3}{|c|}{$\boldsymbol{2(x^6 - 2x^4 + 4x^3 - 2x^2 + 1)}$}\\
\hline
no. & quiver & relations\\
\hline
$11$ & \begin{tabular}{c} \\ \includegraphics[scale=0.8]{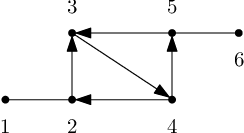} \\ \\ \end{tabular} & \begin{tabular}{c} \\
$(2,3,4)$, $(3,4,2)$, $(5,3,4)$, $(3,4,5)$,\\ $(4,2,3)-(4,5,3)$ \\ \end{tabular}\\
\hline
\end{tabular}
\end{center}}

{\scriptsize
\begin{longtable}[t]{|c|p{5cm}<{\centering}|p{5cm}<{\centering}|}
\hline
\multicolumn{3}{|c|}{}\\[-6pt]
\multicolumn{3}{|c|}{$\boldsymbol{2(x^6 - x^4 + 2x^3 - x^2 + 1)}$}\\
\hline
no. & quiver & relations\\
\hline
\endfirsthead
\hline
no. & quiver & relations\\
\hline
\endhead
$2$ & \begin{tabular}{c} \\ \includegraphics[scale=0.8]{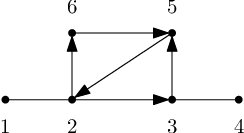} \\ \\ \end{tabular} & \begin{tabular}{c} \\
$(3,5,2)$, $(5,2,3)$, $(6,5,2)$, $(5,2,6)$, \\ $(2,6,5)-(2,3,5)$ \\ \end{tabular}\\
\hline
$7$ & \begin{tabular}{c} \\ \includegraphics[scale=0.8]{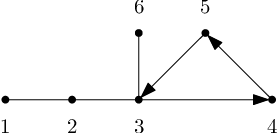} \\ \\ \end{tabular} & \begin{tabular}{c} \\
$(3,4,5)$, $(4,5,3)$, $(5,3,4)$  \\ \\ \end{tabular}\\
\hline
$12$ & \begin{tabular}{c} \\ \includegraphics[scale=0.8]{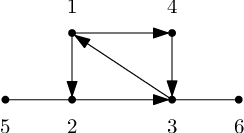} \\ \\ \end{tabular} & \begin{tabular}{c} \\
$(2,3,1)$, $(3,1,2)$, $(4,3,1)$, $(3,1,4)$, \\ $(1,2,3)-(1,4,3)$ \\ \end{tabular}\\
\hline
\end{longtable}}

{\scriptsize
\begin{longtable}[t]{|c|p{5cm}<{\centering}|p{5cm}<{\centering}|}
\hline
\multicolumn{3}{|c|}{}\\[-6pt]
\multicolumn{3}{|c|}{$\boldsymbol{3(x^6 + x^3 + 1)}$}\\
\hline
no. & quiver & relations\\
\hline
\endfirsthead
\hline
no. & quiver & relations\\
\hline
\endhead
$3$ & \begin{tabular}{c} \\ \includegraphics[scale=0.8]{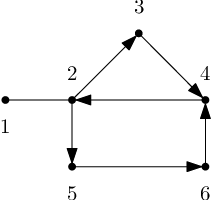} \\ \\ \end{tabular} & \begin{tabular}{c} \\
$(3,4,2)$, $(4,2,3)$, $(5,6,4,2)$, $(6,4,2,5)$, \\ $(4,2,5,6)$, $(2,3,4)-(2,5,6,4)$ \\ \end{tabular}\\
\hline
$4$ & \begin{tabular}{c} \\ \includegraphics[scale=0.8]{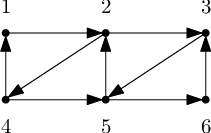} \\ \\ \end{tabular} & \begin{tabular}{c}
$(1,2,4)$, $(2,4,1)$, $(5,2,4)$, $(3,5,2)$, \\ $(3,5,6)$, $(6,3,5)$, $(2,4,5)-(2,3,5)$, \\
$(4,1,2)-(4,5,2)$, $(5,2,3)-(5,6,3)$ \\ \end{tabular}\\
\hline
$6$ & \begin{tabular}{c} \\ \includegraphics[scale=0.8]{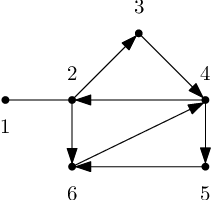}  \\ \\ \end{tabular}& \begin{tabular}{c}
$(3,4,2)$, $(4,2,3)$, $(6,4,2)$, $(5,6,4)$, \\ $(6,4,5)$, $(2,3,4)-(2,6,4)$, \\ $(4,2,6)-(4,5,6)$ \\ \end{tabular}\\
\hline
$10$ & \begin{tabular}{c} \\ \includegraphics[scale=0.8]{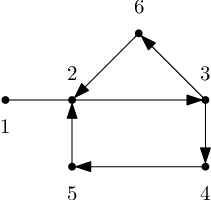} \\ \\ \end{tabular} & \begin{tabular}{c} \\
$(2,3,6)$, $(6,2,3)$, $(4,5,2,3)$, $(5,2,3,4)$, \\ $(2,3,4,5)$, $(3,6,2)-(3,4,5,2)$ \\ \end{tabular}\\
\hline
$14$ & \begin{tabular}{c} \\ \includegraphics[scale=0.8]{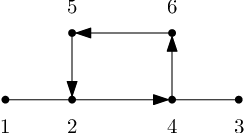} \\ \\ \end{tabular} & \begin{tabular}{c} \\
$(2,4,6,5)$, $(4,6,5,2)$, $(6,5,2,4)$, \\ $(5,2,4,6)$ \\ \end{tabular}\\
\hline
$16$ & \begin{tabular}{c} \\ \includegraphics[scale=0.8]{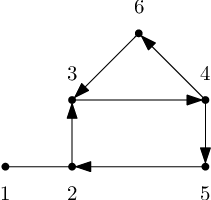} \\ \\ \end{tabular} & \begin{tabular}{c} \\
$(3,4,6)$, $(6,3,4)$, $(5,2,3,4)$, $(2,3,4,5)$, \\ $(3,4,5,2)$, $(4,6,3)-(4,5,2,3)$ \\ \end{tabular}\\
\hline
$19$ & \begin{tabular}{c} \\ \includegraphics[scale=0.8]{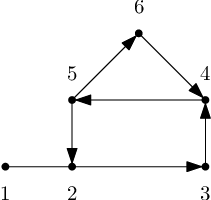} \\ \\ \end{tabular} & \begin{tabular}{c} \\
$(4,5,6)$, $(6,4,5)$, $(2,3,4,5)$, $(3,4,5,2)$, \\ $(4,5,2,3)$, $(5,6,4)-(5,2,3,4)$ \\ \end{tabular}\\
\hline
$20$ & \begin{tabular}{c} \\ \includegraphics[scale=0.8]{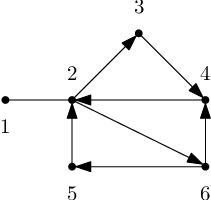} \\ \\ \end{tabular} & \begin{tabular}{c}
$(3,4,2)$, $(4,2,3)$, $(4,2,6)$, $(2,6,5)$, \\ $(5,2,6)$, $(2,3,4)-(2,6,4)$, \\ $(6,4,2)-(6,5,2)$ \\ \end{tabular}\\
\hline
$21$ & \begin{tabular}{c} \\ \includegraphics[scale=0.8]{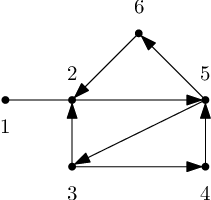} \\ \\ \end{tabular} & \begin{tabular}{c}
$(2,5,6)$, $(6,2,5)$, $(2,5,3)$, $(4,5,3)$, \\ $(5,3,4)$, $(5,6,2)-(5,3,2)$, \\ $(3,2,5)-(3,4,5)$ \\ \end{tabular}\\
\hline
\end{longtable}}

{\scriptsize
\begin{longtable}[t]{|c|p{5cm}<{\centering}|p{5cm}<{\centering}|}
\hline
\multicolumn{3}{|c|}{}\\[-6pt]
\multicolumn{3}{|c|}{$\boldsymbol{4(x^6 + x^4 + x^2 + 1)}$}\\
\hline
no. & quiver & relations\\
\hline
\endfirsthead
\hline
no. & quiver & relations\\
\hline
\endhead
$5$ & \begin{tabular}{c} \\ \includegraphics[scale=0.8]{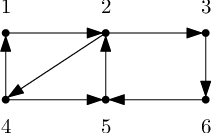}  \\ \\ \end{tabular}& \begin{tabular}{c}
$(1,2,4)$, $(2,4,1)$, $(5,2,4)$, $(3,6,5,2)$, \\ $(6,5,2,3)$, $(5,2,3,6)$, $(4,1,2)-(4,5,2)$, \\ $(2,4,5)-(2,3,6,5)$ \\ \end{tabular}\\
\hline
$9$ & \begin{tabular}{c} \\ \includegraphics[scale=0.8]{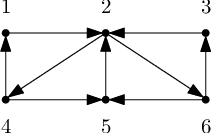} \\ \\ \end{tabular} & \begin{tabular}{c}
$(1,2,4)$, $(2,4,1)$, $(5,2,4)$, $(5,2,6)$, \\ $(2,6,3)$, $(3,2,6)$, $(4,1,2)-(4,5,2)$, \\ $(2,4,5)-(2,6,5)$, $(6,3,2)-(6,5,2)$ \\ \end{tabular}\\
\hline
$13$ & \begin{tabular}{c} \\ \includegraphics[scale=0.9]{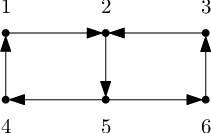} \\ \\ \end{tabular} & \begin{tabular}{c}
$(1,2,5,4)$, $(2,5,4,1)$, $(4,1,2,5)$, \\ $(3,2,5,6)$, $(2,5,6,3)$, $(6,3,2,5)$, \\ $(5,4,1,2)-(5,6,3,2)$ \\ \end{tabular}\\
\hline
$15$ & \begin{tabular}{c} \\ \includegraphics[scale=0.8]{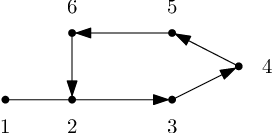} \\ \\ \end{tabular} & \begin{tabular}{c} \\
$(2,3,4,5,6)$, $(3,4,5,6,2)$, $(4,5,6,2,3)$, \\ $(5,6,2,3,4)$, $(6,2,3,4,5)$ \\ \end{tabular}\\
\hline
$17$ & \begin{tabular}{c} \\ \includegraphics[scale=0.8]{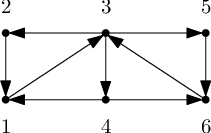} \\ \\ \end{tabular} & \begin{tabular}{c}
$(1,3,2)$, $(2,1,3)$, $(1,3,4)$, $(6,3,4)$, \\ $(6,3,5)$, $(5,6,3)$, $(3,2,1)-(3,4,1)$, \\ $(3,4,6)-(3,5,6)$, $(4,1,3)-(4,6,3)$ \\ \end{tabular}\\
\hline
$18$ & \begin{tabular}{c} \\ \includegraphics[scale=0.8]{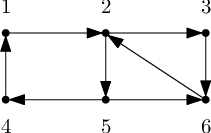} \\ \\ \end{tabular} & \begin{tabular}{c}
$(3,6,2)$, $(6,2,3)$, $(6,2,5)$, $(1,2,5,4)$, \\ $(2,5,4,1)$, $(4,1,2,5)$, $(2,3,6)-(2,5,6)$, \\ $(5,6,2)-(5,4,1,2)$ \\ \end{tabular}\\
\hline
\end{longtable}}

{\scriptsize
\begin{center}
\begin{tabular}[t]{|c|p{5cm}<{\centering}|p{5cm}<{\centering}|}
\hline
\multicolumn{3}{|c|}{}\\[-6pt]
\multicolumn{3}{|c|}{$\boldsymbol{4(x^6 + x^5 - x^4 + 2x^3 - x^2 + x + 1)}$}\\
\hline
no. & quiver & relations\\
\hline
$8$ & \begin{tabular}{c} \\ \includegraphics[scale=0.8]{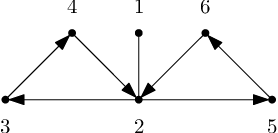}  \\ \\ \end{tabular}& \begin{tabular}{c} \\
$(2,3,4)$, $(3,4,2)$, $(4,2,3)$, \\ $(2,5,6)$, $(5,6,2)$, $(6,2,5)$ \\ \end{tabular}\\
\hline
\end{tabular}
\end{center}}

The rest of this section is devoted to proving our main theorem (as
stated in the introduction) for type $E_6$. In order to do this, it is
enough to show the equivalence of conditions~(a) and~(d) in that
theorem, namely, that any two cluster-tilted algebras of type $E_6$
with the same associated polynomial are connected by a sequence of good
mutations.

First, we give a detailed example. Then we list the other good
mutations in a shorter form which is explained after the example.

\begin{eg}
Let $A_7$ be the cluster-tilted algebra corresponding to one
representative of the quiver number $7$ (as depicted in the left
picture). After mutation at vertex $3$ we get a representative of the
quiver number $2$ as in the right picture. We denote the corresponding
cluster-tilted algebra by $A_2$.

\begin{center}
\begin{tabular}{cc}
\includegraphics[scale=1.2]{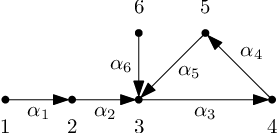} & \qquad \includegraphics[scale=1.2]{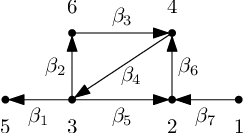}\\
\end{tabular}
\end{center}

Invoking the algorithm in~\cite{BMR_finite}, we see that the relations
of $A_7$ are given by the three zero-relations $\alpha_3\alpha_4$,
$\alpha_4\alpha_5$ and $\alpha_5\alpha_3$. Similarly, the
cluster-tilted algebra $A_2$ has the zero-relations
$\beta_3\beta_4$, $\beta_4\beta_2$, $\beta_4\beta_5$, $\beta_6\beta_4$
and the commutativity-relation $\beta_2\beta_3 = \beta_5\beta_6$. The
corresponding Cartan matrices are therefore

\begin{center}
\begin{tabular}{ccc}
$C_{A_7} = \left(\begin{array}[c]{cccccc} 1&1&1&1&0&0\\ 0&1&1&1&0&0\\
0&0&1&1&0&0\\ 0&0&0&1&1&0\\ 0&0&1&0&1&0\\ 0&0&1&1&0&1  \end{array}\right)$ &
and &
$C_{A_2} = \left(\begin{array}[c]{cccccc} 1&1&0&1&0&0\\ 0&1&0&1&0&0\\
0&1&1&1&1&1\\ 0&0&1&1&1&0\\ 0&0&0&0&1&0\\ 0&0&0&1&0&1  \end{array}\right)$ ,
\end{tabular}
\end{center}
and we compute the associated polynomial to be $2(x^6 - x^4 + 2x^3 -
x^2 +1)$.

We examine now the algebra mutations at the vertex $3$. Since the arrow
$\alpha_6 : 6 \to 3$ does not appear in any relation of $A_7$, its
composition with any non-zero path starting at vertex $3$ is non-zero.
Thus the negative mutation $\mu_3^-(A_7)$ is defined. Since the arrow
$\beta_1 : 3 \to 5$ does not appear in any relation of $A_2$, its
composition with any non-zero path ending at vertex $3$ is non-zero.
Thus the positive mutation $\mu_3^+(A_2)$ is defined.

Since the quiver of $A_2$ is obtained by mutating the quiver of $A_7$
at the vertex $3$, condition (a) of Proposition~\ref{p:goodmut} holds,
hence the quiver mutation at vertex $3$ is good and the corresponding
cluster-tilted algebras $A_7$ and $A_2$ are derived equivalent.
\end{eg}

We now write enough good mutations so that the reader can easily verify
that the quivers of any two cluster-tilted algebras with the same
associated polynomial are indeed connected by a sequence of these
mutations.

We present these good mutations in a concise form in the tables below.
The algebras are numbered according
to the numbering of their quivers, i.e. $A_i$ denotes the
cluster tilted algebra corresponding to the quiver no.\,$i$ in the above tables.
We list the
cluster-tilted algebra together with an orientation of the arrows
incident to the mutated vertex (if they can be oriented arbitrarily),
the vertex we mutate at and the resulting cluster-tilted algebra. Note
that up to a permutation of the vertices and sink/source equivalence,
it is one of the $21$ in our list, and the necessary permutations (written
as a product of disjoint cycles) are given in each case.
By Lemma~\ref{l:goodmutop}, the
mutation of the opposite quivers at the same vertex is also good, so we
also indicate this under ``opposite case'' when the corresponding
algebras are not isomorphic to the original ones. The verification that
all these mutations are good is done as in the previous example.
\begin{center}
\begin{tabular}{|cccc|}
\hline
\multicolumn{4}{|c|}{}\\[-10pt]
\multicolumn{4}{|c|}{$\boldsymbol{2(x^6 - x^4 + 2x^3 -x^2 + 1)}$}\\
\hline
algebra (orientation) & vertex & mutated algebra & permutation\\
\hline
$A_7$ ($2 \to 3$, $6 \to 3$) & $3$ & $A_2$ & $(145)(23)$\\
$A_2$ ($1 \to 2$) & $2$ & $A_{12}$ & $(23)(456)$\\
\hline
\end{tabular}
\end{center}

\begin{center}
\begin{tabular}{|ccccc|}
\hline
\multicolumn{5}{|c|}{}\\[-10pt]
\multicolumn{5}{|c|}{$\boldsymbol{3(x^6 + x^3 + 1)}$}\\
\hline
algebra (orientation) & vertex & mutated algebra & permutation & opposite case\\
\hline
$A_3$ & $5$ & $A_{20}$ & $\mathrm{id}$ & $A_{10}$ $\sim$ $A_{20}$\\
$A_3$ ($1 \to 2$) & $2$ & $A_4$ & $(143)$ & $A_{10}$ $\sim$ $A_4$\\
$A_{20}$ & $4$ & $A_{14}$ & $\mathrm{id}$ & \\
$A_6$ & $5$ & $A_3$ & $(56)$ & $A_{21}$ $\sim$ $A_{10}$\\
$A_{16}$ & $5$ & $A_6$ & $(365)$ & $A_{19}$ $\sim$ $A_{21}$\\
\hline
\end{tabular}
\end{center}

\begin{center}
\begin{tabular}{|ccccc|}
\hline
\multicolumn{5}{|c|}{}\\[-10pt]
\multicolumn{5}{|c|}{$\boldsymbol{4(x^6 + x^4 + x^2 + 1)}$}\\
\hline
algebra (orientation) & vertex & mutated algebra & permutation & opposite case\\
\hline
$A_5$ & $3$ & $A_9$ & $\mathrm{id}$ & $A_{18}$ $\sim$ $A_{17}$\\
$A_{15}$ ($1 \to 2$) & $2$ & $A_5$ & $(243)(56)$ & $A_{15}$ $\sim$ $A_{18}$\\
$A_{13}$ & $4$ & $A_5$ & $(14)(25)(36)$ &\\
\hline
\end{tabular}
\end{center}


\section{Standard forms for the derived equivalence classes}

In this section we provide standard forms for the derived equivalence
classes of cluster-tilted algebras of Dynkin type $E$.
For type $E_6$ the derived equivalence classification has been carried out
in Section \ref{Sec:E6}. For more details
on the actual classification in types $E_7$ and $E_8$ we refer the
reader to the supplementary material \cite{Sefi:www}
as well as to earlier versions of the paper~\cite{BHLTypeE} on the arXiv.

Given such derived equivalence class, we calculate its standard form in
the following way. We start by considering the cluster-tilted algebras
within that class whose quivers have minimal number of arrows. It turns
out that except for one derived class in type $E_7$ consisting of a
single algebra, the quivers with relations of all these ``minimal''
algebras can be described as iterated gluing of quivers with relations
of three possible kinds:
\begin{enumerate}
\item \label{en:cycle}
\emph{Cycles}, consisting of $n$ vertices $a_0, a_1, \dots, a_{n-1}$
and $n$ arrows $a_i \to a_{i+1}$ (where indices are taken modulo $n$).
The relations are given by the $n$ paths $a_i, a_{i+1}, \dots,
a_{i+n-1}$ of length $n-1$ (again indices are taken modulo $n$) for $0
\leq i < n$. Note that these are actually cluster-tilted algebras of
type $D_n$~\cite{BHL10,Vatne}. The values of $n$ which occur are
$3,4,5,6,7$.

\item \label{en:tripenta}
The two cluster-tilted algebras of types $D_4$ and $D_6$ given by the
quivers with relations
\begin{align*}
&\begin{array}{c}
\xymatrix@=0.3pc{
&& {\bullet_1} \ar[dll] \ar[drr] \\
{\bullet_2} \ar[drr] && && {\bullet_3} \ar[dll] \\
&& {\bullet_4} \ar[uu]
}
\end{array}
& &
\begin{array}{l}
(1,2,4)-(1,3,4), (2,4,1), (3,4,1), (4,1,2), (4,1,3)
\end{array}
\intertext{and}
&\begin{array}{c}
\xymatrix@=0.3pc{
&& && {\bullet_3} \ar[ddr] \\
&& {\bullet_1} \ar[dll] \ar[urr] \\
{\bullet_2} \ar[drr] && &&& {\bullet_4} \ar[ddl] \\
&& {\bullet_6} \ar[uu] \\
&& && {\bullet_5} \ar[ull]
}
\end{array}
& &
\begin{array}{l}
(1,2,6)-(1,3,4,5,6), (2,6,1), (6,1,2), \\
(3,4,5,6,1), (4,5,6,1,3), (5,6,1,3,4), (6,1,3,4,5).
\end{array}
\end{align*}

\item
\emph{Tails}, which are orientations of a Dynkin diagram of type $A_n$
with $n \geq 2$.
\end{enumerate}

Recall that a \emph{gluing} of two quivers with relations $Q_1$ and
$Q_2$ with respect to some chosen vertices $v_1$ in $Q_1$ and $v_2$ in
$Q_2$ is obtained by taking the disjoint union $Q_1 \sqcup Q_2$ and
identifying the vertices $v_1$ and $v_2$. The relations are induced
from those of $Q_1$ and $Q_2$ (i.e.\ there are no additional
relations).

We can consider the minimal algebras up to sink/source equivalence and
taking opposites, as these operations do not change the derived
equivalence class (cf.\ Lemma~\ref{l:goodmutss} and
Corollary~\ref{c:opder}). For many derived equivalence classes this determines
uniquely a standard form within the class. For the classes where this
is not the case, imposing the additional condition that the quiver
should have minimal number of ``free'' tails (i.e.\ tails where no
quiver of the above kinds~\ref{en:cycle} and~\ref{en:tripenta} is glued
to one of their ends) determines uniquely a standard form.

In the tables below we provide lists of these standard forms along with
their associated polynomials. Arrows which are displayed without
orientation can be oriented arbitrarily, without changing the derived
equivalence class.


\bigskip

\begin{longtable}{|c|r||c|r|}
\hline
\multicolumn{4}{|c|}{Representatives of derived equivalence classes for type $E_6$} \\
\hline
\multicolumn{1}{|c|}{Quiver and associated polynomial} & \multicolumn{1}{c||}{\#}  &
\multicolumn{1}{c|}{Quiver and associated polynomial} & \multicolumn{1}{c|}{\#} \\
\hline
\endfirsthead
\hline
\multicolumn{1}{|c|}{Quiver and associated polynomial} & \multicolumn{1}{c||}{\#}  &
\multicolumn{1}{c|}{Quiver and associated polynomial} & \multicolumn{1}{c|}{\#} \\
\hline
\endhead
\begin{tabular}{c}
\\
\includegraphics[scale=0.75]{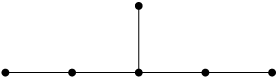} \\ ~ \\
$x^6-x^5+x^3-x+1$ \end{tabular}
& 20 &
\begin{tabular}{c}
\\
\includegraphics[scale=0.75]{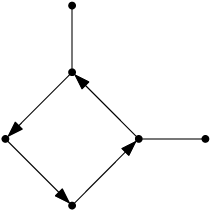} \\ ~ \\
$3(x^6+x^3+1)$ \end{tabular} & 19 \\
\hline
\begin{tabular}{c}
\\
\includegraphics[scale=0.75]{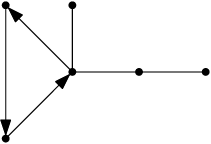} \\ ~ \\
$2(x^6-x^4+2x^3-x^2+1)$ \end{tabular}
& 16 &
\begin{tabular}{c}
\\
\includegraphics[scale=0.75]{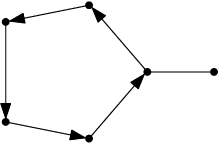} \\ ~ \\
$4(x^6+x^4+x^2+1)$ \end{tabular} & 7 \\
\hline
\begin{tabular}{c}
\\
\includegraphics[scale=0.75]{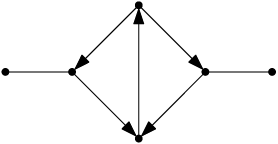} \\ ~ \\
$2(x^6-2x^4+4x^3-2x^2+1)$ \end{tabular}
& 3 &
\begin{tabular}{c}
\\
\includegraphics[scale=0.75]{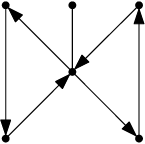} \\ ~ \\
$4(x^6+x^5-x^4+2x^3-x^2+x+1)$ \end{tabular} & 2 \\
\hline
\end{longtable}

\begin{longtable}{|c|r||c|r|}
\hline
\multicolumn{4}{|c|}{Representatives of derived equivalence classes for type $E_7$} \\
\hline
\multicolumn{1}{|c|}{Quiver and associated polynomial} & \multicolumn{1}{c||}{\#}  &
\multicolumn{1}{c|}{Quiver and associated polynomial} & \multicolumn{1}{c|}{\#} \\
\hline
\endfirsthead
\hline
\multicolumn{1}{|c|}{Quiver and associated polynomial} & \multicolumn{1}{c||}{\#}  &
\multicolumn{1}{c|}{Quiver and associated polynomial} & \multicolumn{1}{c|}{\#} \\
\hline
\endhead
\begin{tabular}{c}
\\
\includegraphics[scale=0.75]{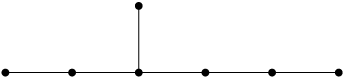} \\ ~ \\
$x^7-x^6+x^4-x^3+x-1$ \end{tabular}
& 64 &
\begin{tabular}{c}
\\
\includegraphics[scale=0.75]{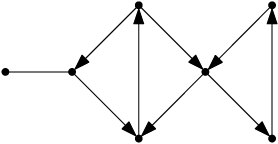} \\ ~ \\
$4(x^7+x^6-2x^5+2x^4-2x^3+2x^2-x-1)$ \end{tabular} & 2 \\
\hline
\begin{tabular}{c}
\\
\includegraphics[scale=0.75]{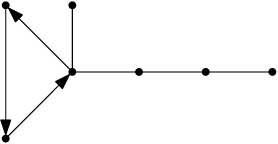} \\ ~ \\
$2(x^7-x^5+2x^4-2x^3+x^2-1)$ \end{tabular}
& 32 &
\begin{tabular}{c}
\\
\includegraphics[scale=0.75]{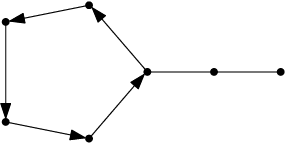} \\ ~ \\
$4(x^7+x^5-x^4+x^3-x^2-1)$ \end{tabular} & 56 \\
\hline
\begin{tabular}{c}
\\
\includegraphics[scale=0.75]{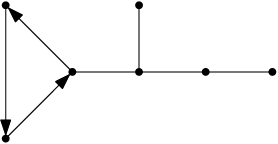} \\ ~ \\
$2(x^7-x^5+x^4-x^3+x^2-1)$ \end{tabular}
& 72 &
\begin{tabular}{c}
\\
\includegraphics[scale=0.75]{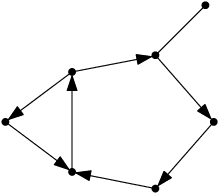} \\ ~ \\
$4(x^7+x^5-2x^4+2x^3-x^2-1)$ \end{tabular} & 8 \\
\hline
\begin{tabular}{c}
\\
\includegraphics[scale=0.75]{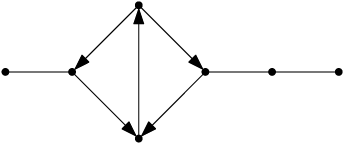} \\ ~ \\
$2(x^7-2x^5+4x^4-4x^3+2x^2-1)$ \end{tabular}
& 8 &
\begin{tabular}{c}
\\
\includegraphics[scale=0.75]{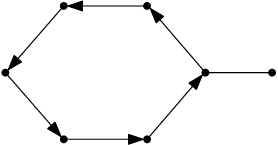} \\ ~ \\
$5(x^7+x^5-x^4+x^3-x^2-1)$ \end{tabular} & 17 \\
\hline
\begin{tabular}{c}
\\
\includegraphics[scale=0.75]{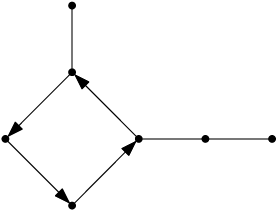} \\ ~ \\
$3(x^7-1)$ \end{tabular}
& 124 &
\begin{tabular}{c}
\\
\includegraphics[scale=0.75]{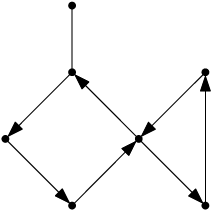} \\ ~ \\
$6(x^7+x^6-x^4+x^3-x-1)$ \end{tabular} & 11 \\
\hline
\begin{tabular}{c}
\\
\includegraphics[scale=0.75]{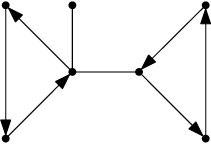} \\ ~ \\
$4(x^7+x^6-x^5+x^4-x^3+x^2-x-1)$ \end{tabular}
& 16 &
\begin{tabular}{c}
\\
\includegraphics[scale=0.75]{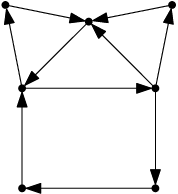} \\ ~ \\
$6(x^7+x^5-x^2-1)$ \end{tabular} & 1 \\
\hline
\begin{tabular}{c}
\\
\includegraphics[scale=0.75]{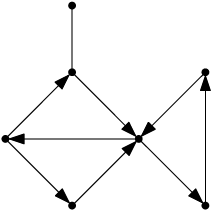} \\ ~ \\
$4(x^7+x^6-x^5-x^4+x^3+x^2-x-1)$ \end{tabular}
& 4 &
\begin{tabular}{c}
\\
\includegraphics[scale=0.75]{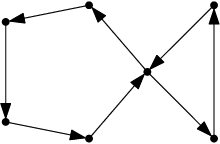} \\ ~ \\
$8(x^7+x^6+x^5-x^4+x^3-x^2-x-1)$ \end{tabular} & 1   \\
\hline
\end{longtable}

\begin{longtable}{|c|r||c|r|}
\hline
\multicolumn{4}{|c|}{Representatives of derived equivalence classes for type $E_8$} \\
\hline
\multicolumn{1}{|c|}{Quiver and associated polynomial} & \multicolumn{1}{c||}{\#}  &
\multicolumn{1}{c|}{Quiver and associated polynomial} & \multicolumn{1}{c|}{\#} \\
\hline
\endfirsthead
\hline
\multicolumn{1}{|c|}{Quiver and associated polynomial} & \multicolumn{1}{c||}{\#}  &
\multicolumn{1}{c|}{Quiver and associated polynomial} & \multicolumn{1}{c|}{\#} \\
\hline
\endhead
\begin{tabular}{c}
\\
\includegraphics[scale=0.75]{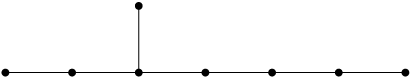} \\ ~ \\
$x^8-x^7+x^5-x^4+x^3-x+1$ \end{tabular}
& 128 &
\begin{tabular}{c}
\\
\includegraphics[scale=0.75]{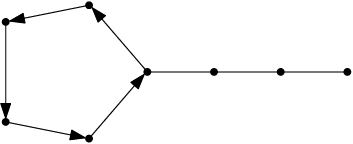} \\ ~ \\
$4(x^8+x^6-x^5+2x^4-x^3+x^2+1)$ \end{tabular} & 221 \\
\hline
\begin{tabular}{c}
\\
\includegraphics[scale=0.75]{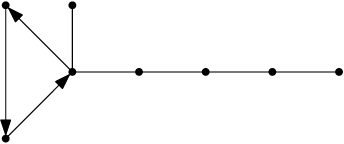} \\ ~ \\
$2(x^8-x^6+2x^5-2x^4+2x^3-x^2+1)$ \end{tabular}
& 64 &
\begin{tabular}{c}
\\
\includegraphics[scale=0.75]{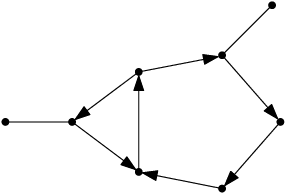} \\ ~ \\
$4(x^8+x^6-2x^5+4x^4-2x^3+x^2+1)$ \end{tabular} & 22 \\
\hline
\begin{tabular}{c}
\\
\includegraphics[scale=0.75]{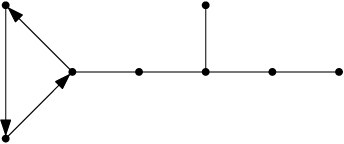} \\ ~ \\
$2(x^8-x^6+x^5+x^3-x^2+1)$ \end{tabular}
& 256 &
\begin{tabular}{c}
\\
\includegraphics[scale=0.75]{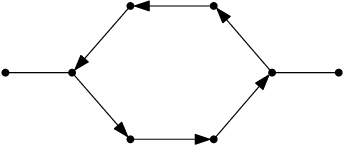} \\ ~ \\
$5(x^8+x^6+x^4+x^2+1)$ \end{tabular} & 167 \\
\hline
\begin{tabular}{c}
\\
\includegraphics[scale=0.75]{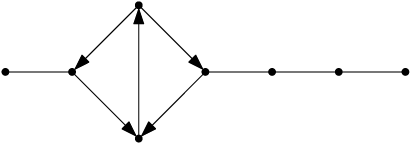} \\ ~ \\
$2(x^8-2x^6+4x^5-4x^4+4x^3-2x^2+1)$ \end{tabular}
& 16 &
\begin{tabular}{c}
\\
\includegraphics[scale=0.75]{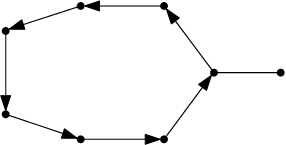} \\ ~ \\
$6(x^8+x^6+x^5+x^3+x^2+1)$ \end{tabular} & 38 \\
\hline
\begin{tabular}{c}
\\
\includegraphics[scale=0.75]{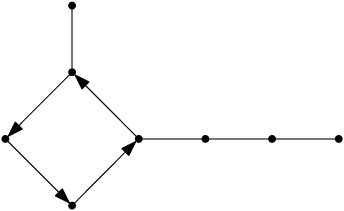} \\ ~ \\
$3(x^8+x^4+1)$ \end{tabular}
& 384 &
\begin{tabular}{c}
\\
\includegraphics[scale=0.75]{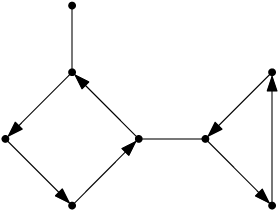} \\ ~ \\
$6(x^8+x^7+2x^4+x+1)$ \end{tabular} & 118 \\
\hline
\begin{tabular}{c}
\\
\includegraphics[scale=0.75]{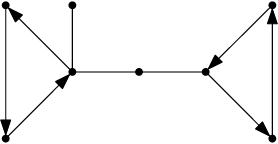} \\ ~ \\
$4(x^8+x^7-x^6+x^5+x^3-x^2+x+1)$ \end{tabular}
& 72 &
\begin{tabular}{c}
\\
\includegraphics[scale=0.75]{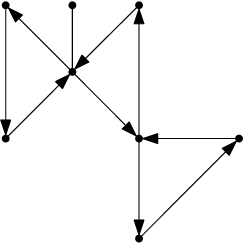} \\ ~ \\
$8(x^8+2x^7+2x^4+2x+1)$ \end{tabular} & 4 \\
\hline
\begin{tabular}{c}
\\
\includegraphics[scale=0.75]{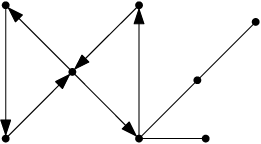} \\ ~ \\
$4(x^8+x^7-x^6+2x^4-x^2+x+1)$ \end{tabular}
& 48 &
\begin{tabular}{c}
\\
\includegraphics[scale=0.75]{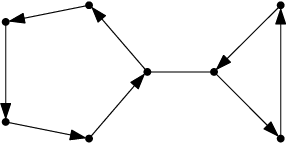} \\ ~ \\
$8(x^8+x^7+x^6+2x^4+x^2+x+1)$ \end{tabular} & 24   \\
\hline
\begin{tabular}{c}
\\
\includegraphics[scale=0.75]{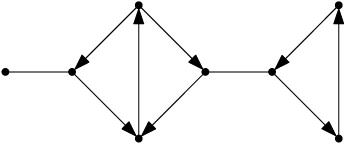} \\ ~ \\
$4(x^8+x^7-2x^6+2x^5+2x^3-2x^2+x+1)$ \end{tabular} & 12 & & \\
\hline
\end{longtable}


\end{document}